\documentclass[12pt]{article}
\usepackage{amssymb,amsbsy,amsmath,amsfonts,amssymb}
\usepackage{latexsym,euscript,exscale}

\usepackage{times}

\newcommand{\pf}{

\smallskip

\noindent {\it Proof : }}

\newcommand {\N}{\mathbb N}

\newcommand {\R}{\mathbb R}

\newcommand {\Z}{\mathbb Z}

\newcommand {\C}{\mathbb C}
\newcommand {\HH}{\mathbb H}

\newcommand {\DD}{{\mathbb D}_2}

\newcommand{\pff}{$\hfill \square$
\smallskip}

\newcommand{\norm}[1]{\ensuremath{\left\|#1\right\|}}

\newtheorem{prop}{Proposition}
\newtheorem{lemm}[prop]{Lemma}
\newtheorem{theo}[prop]{Theorem}

\newtheorem{coro}[prop]{Corollary}
\newtheorem{rema}[prop]{Remark}

\title{Real hereditarily indecomposable Banach spaces and uniqueness of complex
structure}

\author{Valentin Ferenczi}

\date{ }

\begin{document}

\maketitle

\begin{abstract}
There exists a real hereditarily indecomposable Banach space $X=X(\C)$ (resp. $X=X(\HH)$)
 such that
the algebra ${\cal L}(X)/{\cal S}(X)$ is isomorphic to $\C$ (resp. to the
quaternionic division algebra $\HH$).

  Up to isomorphism, $X(\C)$ has exactly 
two complex structures, which are conjugate, totally incomparable,
 and both hereditarily indecomposable. So there exist two Banach spaces which
are isometric as real spaces but totally incomparable as complex spaces.
This extends  results of J. Bourgain and S. Szarek \cite{B,S,S2}. 

The quaternionic example $X(\HH)$, on the other hand, has unique complex structure up to isomorphism; 
there also exists
a space with an unconditional basis, non isomorphic to $l_2$, with unique complex structure. This
answers a question of S. Szarek in \cite{S2}. 
\end{abstract}

\section{Introduction}

It is  well-known  that  any two real Banach spaces which
are isometric must be linearly isometric: this was proved in 1932 
by S. Mazur and S. Ulam \cite{MU}. Recently, G. Godefroy and
N.J. Kalton also proved that if a separable real Banach space embeds isometrically
into a Banach space, then
it embeds linearly isometrically into it \cite{GK}. 

In 1986, J. Bourgain and S. Szarek proved that Mazur-Ulam's result is completely false in the
complex case: there exist two Banach spaces which are linearly isometric as real spaces but non
isomorphic as complex spaces \cite{B,S,S2}.
One of the main results of this paper is the following extension of their result. 

\begin{theo} There exist two Banach spaces which are isometric as real spaces, but totally
incomparable as complex spaces. 
\end{theo}

We recall that two spaces
are said to be totally incomparable if no subspace of one if is isomorphic to a subspace of the
other. We shall also show that this result is optimal, in the sense that there does not exist
a family of more than two Banach spaces which are mutually
isomorphic as real spaces but totally incomparable as complex spaces.

\

Our examples are natural modifications of the hereditarily indecomposable
Banach space of W.T. Gowers and B. Maurey \cite{GM1}.
Hereditarily indecomposable (or H.I.) Banach
spaces were discovered in 1991 by these two authors: a space $X$ is H.I. if
no  (closed) subspace of $X$ is decomposable (i.e. can be written as a direct
sum of infinite dimensional subspaces). Equivalently, for any two
subspaces
$Y$,
$Z$ of $X$ and $\epsilon>0$, there exist $y \in Y$, $z \in Z$ such that
$\norm{y}=\norm{z}=1$ and $\norm{y-z}<\epsilon$. Gowers and Maurey gave
the first known example $X_{GM}$ of a H.I. space, both in the real and
the complex case. 
After this result, many other examples of H.I. spaces with various
additional properties were given. They are too numerous to be all cited
here. We refer to
\cite{A} for a list of these examples. Let us mention however the
remarkable result of S. Argyros and A. Tolias
\cite{AT}: for any separable Banach space $X$ not
containing
$\ell_1$, there is a separable H.I. space with a quotient isomorphic to
$X$.

A Banach space $X$ is said to have a Schauder basis $(e_i)_{i \in \N}$ if any vector of $X$ may be written
uniquely as an infinite sum $\sum_{i \in \N}\lambda_i e_i$. The basis $(e_i)$ is unconditional if any permutation
of $(e_i)$ is again a basis. This is equivalent to saying
 that there exists $C<+\infty$ such that for any vector written
$\sum_{i \in \N}\lambda_i e_i$ in $X$, and any scalar sequence $(\mu_i)_{i \in \N}$ such that
$\forall i \in \N, |\mu_i| \leq |\lambda_i|$, the inequality
$\norm{\sum_{i \in \N}\mu_i e_i} \leq C\norm{\sum_{i \in \N}\lambda_i e_i}$
holds.

Classical spaces, such as $c_0$, $l_p$ for $1 \leq p<+\infty$, $L_p$ for  $1<p<+\infty$, and Tsirelson's
space
$T$ have unconditional bases; or contain subspaces with an unconditional basis in the case of $C([0,1])$ or
$L_1$.
 The H.I. spaces of Gowers and Maurey were the
first known examples of spaces not containing any unconditional basic sequence, thus answering an
old open question by the negative. The importance of H.I. spaces also
stems from the famous Gowers' dichotomy theorem
\cite{G1}: any Banach space contains either a subspace with an
unconditional basis or a H.I. subspace. In some sense, H.I. spaces
capture even more of the general structure of all separable Banach spaces
than classical spaces. By
\cite{A2}, any Banach space containing copies of all separable (reflexive)
H.I. spaces must be universal for the class of separable Banach spaces
(i.e. must contain an isomorphic copy of any separable Banach space). On
the other hand, all spaces with an unconditional basis may be embedded
into the unconditional universal space
$U$ of Pe\l czy\'nski (see i.e \cite{LT} about this). The space $U$ has an unconditional
basis and thus is certainly not universal (for example, it does not contain $L_1$ nor a H.I.
subspace).

\

H.I. spaces have interesting properties linked to the space of operators
defined on them. An operator $s \in {\cal L}(Y,Z)$ is strictly singular if
no restriction of $s$ to an infinite dimensional subspace of $Y$ is an
isomorphism into $Z$. Equivalently, for any $\epsilon>0$, any subspace
$Y'$ of $Y$, there exists $y \in Y'$ such that $\norm{s(y)} < \epsilon
\norm{y}$. The ideal of strictly singular operators is denoted ${\cal
S}(Y,Z)$. It is a pertubation ideal for Fredholm operators, we refer to \cite{LT} about this.
Gowers and Maurey proved that any complex H.I. space
$X$ has what we shall call the
{\em $\lambda Id+S$-property}, i.e. 
any operator on $X$ is of the form $\lambda Id+S$, where $\lambda$ is
scalar and $S$ strictly singular. Note however that
spaces with the $\lambda Id+S$-property which are far from being H.I. were also
constructed
\cite{AM}.  It follows from this property that any operator on
$X$ is either strictly singular or Fredholm with index $0$, and so $X$
cannot be isomorphic to a proper subspace (thus the existence of
$X_{GM}$ answers the old hyperplane's problem,
 which had been solved previously by Gowers \cite{G0}).
 Then in \cite{F1} the author
extended the result: if $X$ is complex H.I. and $Y$
is a subspace of
$X$, then every operator from $Y$ into
$X$ is of the form $\lambda i_{Y,X}+s$, where $\lambda$ is scalar,
$i_{Y,X}$ the canonical injection of $Y$ into $X$, and $s \in {\cal
S}(Y,X)$. This property of operators characterizes complex H.I. spaces.

\

When $X$ is real the situation is more involved. From now on, $X_{GM}$ denotes the real
version of the H.I. space of Gowers and Maurey, as
opposed to the complex version $X_{GM}^{\C}$. The real space
$X_{GM}$ has the property that for any
$Y
\subset X_{GM}$,  any operator from $Y$
to
$X_{GM}$ is of the form $\lambda i_{Y,X_{GM}}+s$, where $s$ is strictly
singular \cite{GM1} (note that this property of operators implies the H.I. property).
In general, a real H.I. space
$X$ must satisfy that for all
$Y \subset X$, $\dim {\cal L}(Y,X)/{\cal S}(Y,X) \leq 4$ \cite{F2}. The
converse is false: the space $X=X_{GM} \oplus X_{GM}$ is not H.I. but for
any $Y
\subset X$, $\dim {\cal L}(Y,X)/{\cal S}(Y,X) \leq 4$ (see the proof
in Remark
\ref{remark1}). Also when $X$  is real H.I., the algebra ${\cal
L}(X)/{\cal S}(X)$ is a division algebra isomorphic to $\R$, $\C$ or the
quaternionic division algebra $\HH$ \cite{F2}. This implies easily, by
continuity of the Fredholm index, that any operator on $X$ is either
strictly singular or Fredholm with index $0$ (this was already proved in
\cite{GM1}), and so
$X$ is not isomorphic to a proper subspace.

\

We will show that each of the values  $2$ and $4$ for $\dim{\cal
L}(X)/{\cal S}(X)$ is indeed possible. We shall build versions of
$X_{GM}$ for which the algebra ${\cal L}(X)/{\cal S}(X)$ is isomorphic to
$\C$ or $\HH$. Furthermore, the complex and the quaternionic
examples satisfy $\dim{\cal L}(Y,X)/{\cal S}(Y,X)=2$ and $4$
respectively, for any subspace $Y$ of $X$:

\begin{theo}
 There exists a real H.I. Banach space $X(\C)$ such that for any subspace $Y$ of $X(\C)$,
$\dim{\cal L}(Y,X(\C))/{\cal S}(Y,X(\C))=2$, and such that the algebra ${\cal L}(X(\C))/{\cal S}(X(\C))$
is isomorphic to the complex field $\C$.

 There exists a real H.I. Banach space $X(\HH)$ such that for any subspace $Y$ of $X(\HH)$,
$\dim{\cal L}(Y,X(\HH))/{\cal S}(Y,X(\HH))=4$, and the algebra ${\cal L}(X(\HH))/{\cal S}(X(\HH))$
is isomorphic to $\HH$, the algebra of quaternions.
\end{theo}

\

 The initial idea of this construction was given
to the author by his  Ph. D. advisor B. Maurey in 1996, but the examples
were never checked. Recently the question was asked the author by D.
Kutzarova and S. Argyros and the interest in this subject was revived by the survey
paper of Argyros \cite{A}, see also \cite{AT}.  The author therefore decided to write down the
proof. The complex structure properties of our examples turned out to be quite interesting.

We shall write the  construction of the real H.I. space denoted
$X(\HH)$ in the quaternionic case, assuming familiarity with
Gowers-Maurey type constructions as in
\cite{GM1} and mainly
\cite{GM2}. The reader will adapt the construction for the
example with complex algebra of operators, denoted $X(\C)$. We shall then
proceed to give the proofs of the operators properties in each case.

\

In the following, spaces and subspaces are supposed infinite-dimensional and closed, unless specified otherwise.
If $X$ is a complex Banach space, with scalar
multiplication denoted $(\lambda+i\mu) x$ for $\lambda, \mu \in \R$ and $x \in X$,
its conjugate
$\overline{X}$ is defined as $X$ equipped with the scalar multiplication
$(\lambda+i\mu).x:=(\lambda-i\mu)x$. Note that $X$ and $\overline{X}$ are isometric as
real spaces.

When $X$ is a real Banach space, a {\em complex structure} $X^I$ on $X$ is $X$ seen as a
complex space with the law $$\forall \lambda,\mu \in \R, (\lambda+i\mu).x=(\lambda Id+\mu
I)(x),$$
where $I$ is some operator on $X$ such that $I^2=-Id$, and renormed with the equivalent norm
$$|||x|||=\sup_{\theta \in \R}\norm{\cos\theta x+\sin\theta Ix}.$$
 Note that the conjugate
of
$X^I$ is the complex structure
$X^{-I}$.

It is known that complex structures do not always exist on
a Banach space, even on a uniformly convex one \cite{S2}.  Gowers also
constructed a space with an unconditional basis on which every operator is a strictly singular perturbation
of a diagonal operator \cite{G0,GM2}, and which therefore does not admit complex structure (this answers Pb
7.1. in
\cite{S2}).

Concerning uniqueness, J. Bourgain
\cite{B}, S. Szarek
\cite{S,S2}, used local random techniques to give an example of  a complex Banach space  not isomorphic to its
conjugate; the space is a "gluing" of finite dimensional spaces which are "far" from their
conjugates. 
Therefore there exist spaces which are isometric  as real spaces but not isomorphic as complex
spaces. Later on,
 N. J. 
Kalton
\cite{K} constructed a simple example defined as a twisted Hilbert space. Recently, R. Anisca \cite{An} constructed a
subspace of $L_p, 1 \leq p<2$ which has the same property, and moreover admits continuum
many non-isomorphic complex structures. Note that these examples fail to have an unconditional basis, since
when a complex Banach space $X$ has an unconditional basis
$(e_n)$, the map $c$ defined by 
$c(\sum_{n \in \N}\lambda_n e_n)=\sum_{n \in \N} \overline{\lambda_n} e_n$ is a $\C$-linear isomorphism from $X$ onto
$\overline{X}$.

\

The real spaces
$X(\C)$ and $X(\HH)$ possess
an operator $I$ such that $I^2=-Id$; the associated complex structures are H.I..
 In fact the
space
$X(\C)$  will look a lot like the complex version of Gowers-Maurey's space
$X_{GM}^{\C}$ seen with its $\R$-linear structure.  However
it doesn't seem at all obvious that the properties of $X(\C)$ are
shared by $X_{GM}^{\C}$ seen as a real space.  
Technically, the difference comes from the fact that the functionals used
in the definition of the norm in our construction of
$X(\C)$ are
$\R$-linear but cannot be allowed to be
$\C$-linear when viewed in the complex setting (see Lemma
\ref{hahnbanach} which prevents this).

\

We show that
$X(\C)$, with the complex structure $X^J(\C)$ associated to some canonical operator $J$, is totally
incomparable with its  conjugate. Therefore we have examples of Banach spaces which are
isometric as real spaces but totally incomparable as complex spaces.
Furthermore, it turns out that
$X^J(\C)$ and its conjugate are the only complex structures on $X(\C)$ up to isomorphism.

The space $X(\HH)$, on the other hand,
admits a unique complex structure up to isomorphism. This answers a question of Szarek
from \cite{S}: he asked whether the Hilbert space was the only space with unique complex
structure.

\begin{theo}
 There exists a real H.I. Banach space which admits exactly two complex structures
up to isomorphism. Moreover, these two complex structures are conjugate and totally
incomparable. 

 There also exists a real H.I. Banach space with unique complex structure up to isomorphism.
\end{theo}

\

We shall see that the  space $X(\C)$ is in some sense the only possible example of space
with totally incomparable complex structures: if
a space $X$ admits two totally incomparable complex structures, then these
 structures must be conjugate up to isomorphism and both
saturated with H.I. subspaces.
It  follows that there cannot be more than two mutually totally
incomparable structures on a Banach space. We shall also note that for any $n \in \N$,
the space $X(\C)^n$ admits exactly $n+1$ complex structures up to isomorphism.

\

Our constructions are different from those of the previous authors, al\-though as no\-ted
by  Maurey in \cite{M}, there are
some similarities between the "few operators" properties of Gowers-Maurey's spaces and the
"few operators" properties of the finite-dimensional spaces glued together in Bourgain-Szarek's
example.
 The peculiar complex structure properties of the spaces $X(\C)$ and $X(\HH)$ follow
directly from their few operators properties.
 For example, $J$ and $-J$ are, up
to strictly singular operators, the only operators on
$X(\C)$ whose square is
$-Id$.

\

We also prove that whenever $\{X_i, 1 \leq i \leq N\}$ is a family of pairwise totally incomparable
  real Banach spaces with the $\lambda Id+S$-property, and $n_i, 1 \leq i \leq N$, are integers,
the direct sum
 $\sum_{1 \leq i \leq N} \oplus X_i^{2n_i}$ has a unique complex structure up to isomorphism.
 This provides additional examples
of spaces not isomorphic to a Hilbert space, which have a unique complex structure.
We also provide an unconditional version of the space $X(\C)$:

\begin{theo} There exists a real Banach space with an unconditional basis, not isomorphic to $l_2$, with unique
complex structure up to isomorphism.
\end{theo}

Finally, note that a complex Banach space, which is H.I. as a real space, is always
complex H.I.. 
Indeed if it  contained two $\C$-linear subspaces in
a direct sum, these would in particular form a direct sum of $\R$-linear
subspaces.
We shall show that the converse is false
(the complexification of
the real
$X_{GM}$ will do).

\section{Some properties of H.I. spaces}

The following was proved in \cite{F1}, \cite{F2}. 

\begin{theo}\label{valentin} \cite{F1},\cite{F2}
Let $X$ be a real H.I. space. Then there exists a division algebra $E$
which is isomorphic to either $\R$, $\C$ or $\HH$, and, for $Y \subset X$,
linear embeddings $i_Y$ of $E_Y={\cal L}(Y,X)/{\cal S}(Y,X)$ into $E$.
Let $\leq$ be the relation between subspaces of $X$ defined by $Z \leq
Y$ iff $Z$ embeds into $Y$ by an isomorphism of the form
$i_{ZX}+s$, where $s \in {\cal S}(Z,X)$.
For any $Z \leq
Y$, the canonical restriction map modulo strictly singular operators
$p_{YZ}:E_Y \rightarrow E_Z$ embeds $E_Y$ into $E_Z$ and satisfies
$i_Y=i_Z p_{YZ}$. The algebra $E$ is actually the inductive limit of
the system $(E_Y,p_{YZ})$ under $\leq$, which is a filter relation.
Furthermore the map
$i_X$ embeds
${\cal L}(X)/{\cal S}(X)$ as a subalgebra of
$E$. 
\end{theo}

A technical lemma (Lemma 2 in \cite{F2}) will be very
useful. Given a Banach space $X$, a subspace $Y$ of $X$ is was defined in
\cite{F2} to be {\em quasi-maximal} in
$X$ if $Y+Z$ for $Z$ infinite dimensional in $X$ is never a direct sum.
Equivalently the quotient map from $X$ onto $X/Y$ is strictly singular.
An obvious remark is that a space $X$ is H.I. if and only if any subspace
of
$X$ is quasi-maximal in $X$.

\begin{lemm}\label{quasimaximal} \cite{F1}
let $X$ be a Banach space, let $T$ be an operator from $X$ into some
Banach space and let
$Y$ be quasi-maximal in
$X$. Then $T$ is strictly singular if and only if $T_{|Y}$ is strictly
singular.

In particular, if $X$ is H.I. and $Y$ a subspace of $X$, then
$T$ is strictly singular if and only if $T_{|Y}$ is strictly singular.
\end{lemm}

 The "filter property" will
denote the fact that if $X$ is H.I. and
$Y, Z$ are subspaces of $X$, then there exist a subspace $W$ such that
$W \leq Y$ and $W \leq Z$ - in particular $W$ embeds into $Y$ and $Z$
(Lemma 1 in \cite{F1}). 

We recall that a space $X$ is said to be $HD_n$ if $X$ contains at
most and exactly $n$ infinite dimensional subspaces in a direct sum
\cite{F2}. For example, a space is $HD_1$ if and only if it is H.I..

We finally recall that $X_{GM}$ denotes the real H.I. space of Gowers and
Maurey. The following remark shows that real H.I. spaces are not characterized by
the property that $\dim{\cal L}(Y,X)/{\cal S}(Y,X) \leq 4$ for all subspaces $Y$ of $X$.

\begin{rema}\label{remark1}
Let $X=X_{GM} \oplus X_{GM}$. For
any $Y
\subset X$, $\dim {\cal L}(Y,X)/{\cal S}(Y,X) \leq 4$.
\end{rema}

\pf 
By \cite{F2} Corollary 2, $X$ is $HD_2$ as a direct sum of two H.I. spaces.
Let $Y$ be an arbitrary subspace of $X$, then $Y$ is either H.I. or
$HD_2$. Let $d$ be the dimension of ${\cal L}(Y,X)/{\cal S}(Y,X)$. Assume
$Y$ is
$HD_2$. Then $Y$ contains a direct sum of H.I. spaces $Z_1 \oplus Z_2$. By
\cite{F2} Corollary 3, passing to subspaces and by the
filter property, we may assume that $Z_1$ and $Z_2$ are isomorphic and
embed into $X_{GM}$. Since    $\dim {\cal L}(Z_1,X_{GM})/{\cal
S}(Z_1,X_{GM})=1$, we deduce that
$\dim {\cal L}(Z_1 \oplus Z_2,X)/{\cal S}(Z_1 \oplus Z_2,X) \leq 4$.

Since $Z_1 \oplus Z_2$ is $HD_2$ (\cite{F2} Corollary 2), it is
quasimaximal in $Y$, so the restriction map
$r$ defined from ${\cal L}(Y,X)/{\cal S}(Y,X)$ into 
${\cal L}(Z_1 \oplus Z_2,X)/{\cal S}(Z_1 \oplus Z_2,X)$ by
$r(\tilde{T})=\widetilde{T_{|Z_1 \oplus Z_2}}$ is well defined and
injective (Lemma \ref{quasimaximal}) - here $\tilde{T}$ denotes the class
of $T$ modulo strictly singular operators. It follows that
$$d \leq   
\dim {\cal L}(Z_1 \oplus Z_2,X)/{\cal S}(Z_1 \oplus Z_2,X) = 4.$$
If $Y$ is H.I., we do a similar proof, passing to a subspace $Z$ of $Y$
which embeds into $X_{GM}$, and obtain by the same H.I. properties 
that
$$d \leq  
\dim {\cal L}(Z,X)/{\cal S}(Z,X) =2.$$ \pff

\section{Construction of real H.I. spaces $X(\HH)$ and $X(\C)$}

 Let $c_{00}$ be the vector space of all real sequences which are eventually
$0$. Let $(e_n)_{n \in \N}$ be the standard basis of $c_{00}$.
Given a family of vectors $\{x_i, i \in I\}$, $[x_i, i \in I]$ denotes the closed linear
span of $\{x_i, i \in I\}$. For $k \in \N$, let $F_k=[e_{4k-3},e_{4k-2},e_{4k-1},e_{4k}]$.
The sequence $F_k$ will provide a $4$-dimensional decomposition of
$X(\HH)$.

We proceed to definitions which are adaptations of the usual
Gowers-Maurey definitions to the $4$-dimensional decomposition context. 

If $E \subset \N$, then we shall also use the letter $E$ for the
projection from $c_{00}$ to $c_{00}$ defined by $E(\sum_{i \in
\N}x_i)=\sum_{i \in E}x_i$, where $x_i \in F_i, \forall i \in \N$.
If $E,F \subset \N$, then we write $E<F$ to mean that $\max E<\min F$,
and if $k \in \N$ and $E \subset \N$, then we write $k<E$ to mean
 $k<\min E$. The support of a vector $x=\sum_i x_i \in c_{00}, x_i \in
F_i$, is the set of $i \in \N$ such that $x_i \neq 0$. An interval of
integers is the intersection of an interval of $\R$ with $\N$. The range
of a vector, written $ran(x)$, is the smallest interval containing its
support. We shall write $x<y$ to mean $ran(x)<ran(y)$. If
$x_1<\cdots<x_n$, we shall say that $x_1,\ldots,x_n$ are successive.

The class of functions $\cal F$ is defined as in \cite{GM1}, and the
function
$f \in {\cal F}$ is  defined on
$[1,+\infty)$ by
$f(x)=\log_2(1+x)$.
Let ${\cal X}_4$ be the set of normed spaces of the form
$X=(c_{00},\norm{.})$ such that $(F_i)_{i=1}^{\infty}$ is a monotone
Schauder decomposition of $X$ where each $e_k$ is normalized. If $f \in
\cal F$, $X \in {\cal X}_4$ and every $x \in X$ satisfies the inequality
$$\norm{x} \geq \sup\{f(N)^{-1}\sum_{i=1}^N \norm{E_i x}: N \in \N,
E_1<\cdots<E_n\},$$
where the $E_i$'s are intervals, then we shall say that $X$ 
satisfies a lower $f$-estimate (with respect to $(F_i)_{i=1}^{\infty}$).

Special vectors are considered in
\cite{GM2}. We give their definitions in our context, as well as some
lemmas without proof: indeed their proof is essentially the proof from
\cite{GM2} word by word, with the difference that "successive"
and "lower
$f$-estimate" mean with respect to $(e_i)$ in their case and to $(F_k)$
in ours. As our definitions are also based on the decomposition $(F_k)$
instead of
$(e_i)$, it is easy to check that the proofs are indeed
valid.  We shall only point out the parts of the proofs which require a
non-trivial modification.   Alternatively, in the case of the real H.I.
space with  complex algebra of operators, these lemmas correspond exactly
to the lemmas of
\cite{GM2} for the complex space $X_{GM}^{\C}$, with our
2-dimensional real decomposition interpreted as a Schauder basis on $\C$.

For $X \in {\cal X}_4$, $x \in X$, and every integer $N \geq 1$, we consider the equivalent
norm on $X$ defined by
$$\|x\|_{(N)}= \hbox{sup} \sum_{i=1}^{N} \|E_{i}x\|,$$
where the supremum is over all sequences $E_1, E_2, \cdots, E_N$ 
of successive intervals.

For $0<\epsilon \leq 1$ and $f \in {\cal F}$, we say that a sequence $x_1,\ldots,x_N$ of
suc\-ces\-sive vectors satisifies the RIS($\epsilon$) condition 
(for $f$) if there exists a sequence
$n_1<\ldots<n_N$ of integers such 
that $\norm{x_i}_{(n_i)}
\leq 1$ for each $i=1,\ldots,N$, $n_1>
(2N/f'(1))f^{-1}(N^2/\epsilon^2)$,
and
$\epsilon \sqrt{f(n_i)}>|ran(\sum_{j=1}^{i-1}x_j)|$, for 
$i=2,\ldots,N$. Observe that when $x_1,\ldots,x_N$ satisfies the
RIS($\epsilon$) condition, then $Ex_1,\ldots,Ex_N$ also does for any
interval $E$.

Given $g \in {\cal F}$, $M \in \N$ and $X \in {\cal X}_4$, an $(M,g)$-form on $X$ is a
functional of norm at most $1$ which can be written $\sum_{j=1}^M x_j^*$ for a sequence
$x_1^*<\ldots<x_M^*$ of successive functionals of norm at most $g(M)^{-1}$. Observe that
if $x^*$ is an $(M,g)$-form then $|x^*(x)| \leq g(M)^{-1}\norm{x}_{(M)}$ for any $x$.

\begin{lemm}\label{gm1}
Let $X \in {\cal X}_4$. Let $f,g \in {\cal F}$ be such that $\sqrt{f} \leq
g$. Assume that
$x_1,\ldots,x_N \in X$ satisfies the RIS($\epsilon$)-condition for $f$.
If 
$x^*$ is a $(k,g)$-form for some integer $k \geq 2$ then
$$|x^*(\sum_{i=1}^N x_i)| \leq \epsilon+1+N/\sqrt{f(k)}.$$
In particular, $|x^*(x_1+\ldots+x_N)| < 1+2\epsilon$ when $k>f^{-1}(N^2/\epsilon^2)$.
\end{lemm}

\pf Reproduce the proof of \cite{GM2} Lemma 1,
noting that, for $x \in c_{00}$, $\norm{x}_{c_0}=\max_{i \in \N} {\norm{x_i}}$
 if $x=\sum_{i \in \N} x_i$  with $x_i \in F_i, \forall i \in \N$. \pff

\begin{lemm}\label{gm3} Let $X \in {\cal X}_4$. Let $f,g \in {\cal F}$,
$\sqrt{f}
\leq g$, and let $x_1,\ldots,x_N \in X$ satisfy the RIS($\epsilon$)
condition for
$f$. Let
$x=\sum_{i=1}^N x_i$, and suppose that
$$\norm{Ex} \leq 1 \vee \sup\{|x^*(Ex)|: x^* (k,g)-form, k \geq 2\},$$
for every interval $E$. Then $\norm{x} \leq (1+2\epsilon)Ng(N)^{-1}$.
\end{lemm}

\pf Reproduce the proof of \cite{GM2} Lemma 3, using \cite{GM2} Lemma 4
in its $4$-dimen\-sio\-nal decomposition version. \pff

\begin{lemm}\label{gm4} Let $X \in {\cal X}_4$ satisfy a lower
$f$-estimate. Then for every
$n
\in
\N$ and $\epsilon>0$, every subspace of $X$ generated by a
sequence of successive vectors contains a vector
$x$ of finite support such that $\norm{x}=1$ and $\norm{x}_{(n)} \leq
1+\epsilon$. Hence, for every $N \in \N$, every subspace generated by a
sequence of successive vectors contains a sequence $x_1,\ldots,x_N$
satisfying the RIS($\epsilon$) condition with $\norm{x_i}\geq
(1+\epsilon)^{-1}$.
\end{lemm}

\pf Given a sequence $(x_n)_{n \in \N}$ of successive vectors in $X$,
$(x_n)$ is basic bimonotone, and for vectors in $[x_n, n \in \N]$, the
notions of lower $f$-estimate, successive vectors, etc... with respect to
$(F_k)_{k \in \N}$ correspond  to the usual notions of lower
$f$-estimate, successive vectors, etc... with respect to the basis
$(x_n)$. Therefore the conclusion holds by Lemma 4 in \cite{GM2}. 
\pff

\

We now pass to the definition of $X=X(\HH)$. Let ${\bf Q} \subset
c_{00}$ be the set of sequences with rational coordinates and modulus at
most $1$. Let $J \subset \N$ be a set such that if $m<n$ and
$m,n \in J$, then $\log \log \log n \geq 2m$. We write $J$ in increasing
order as
$\{j_1,j_2,\ldots\}$. We shall also need $f(j_1)>256$ where $f(x)$ is
still the function $\log_2(x+1)$. Let $K,L \subset J$ be the sets
$\{j_1,j_3,\ldots\}$ and $\{j_2,j_4,\ldots\}$.

Let $\sigma$ be an injection from the collection of finite sequences of
successive elements of ${\bf Q}$ to $L$. Given $X \in {\cal X}_4$ and $f
\in {\cal F}$ such that $X$ satisfies a lower $f$-estimate (with respect
to $(F_k)$), and given an integer $m \in \N$, let $A_m^*(X)$ be the set
of functionals of the form $f(m)^{-1}\sum_{i=1}^m x_i^*$ such that
$x_1^*<\ldots<x_m^*$ and $\norm{x_i^*} \leq 1$.

If $k \in \N$, let $\Gamma_k^X$ be the set of sequences
$y_1^*<\cdots<y_k^*$ such that $y_i^* \in {\bf Q}$ for each $i$, $y_1^*
\in A^*_{j_{2k}}(X)$ and
$y^*_{i+1} \in A^*_{\sigma(y^*_1,\ldots,y^*_i)}(X)$ for each $1 \leq i
\leq k-1$. These sequences are called special sequences. If
$(g_i)_{i=1}^k$, $k \in K$, is a special sequence, then the functional
$f(k)^{-1/2}\sum_{j=1}^k g_j$ is a special functional on $X$ of size $k$.
The set of such functionals is denoted $B_k^*(X)$. If $f \in \cal F$ and
$g(k)=f(k)^{1/2}$, then a special functional of size $k$ is also a
$(k,g)$-form.

The quaternionic division algebra may be represented as an algebra of
operators on $\R^4$. It is then generated by a family
$\{Id_{\R^4},u,v,w\}$, where $u,v,w$ satisfy the relations
$u^2=v^2=w^2=-Id_{\R^4}$,
$uv=-vu=w$, $vw=-wv=u$, and $wu=-uw=v$. We may identify $u$, $v$,
and $w$ with operators $u_k,v_k$ and $w_k$ on each $F_k$ using the
identification to
$\R^4$ via the canonical basis $e_{4k-3},\ldots,e_{4k}$ of $F_k$. We then
define linear operators $U,V$ and $W$ on $c_{00}$ by, for all
$k \in \N$, $U_{|F_k}=u_k$ (resp.
$V_{|F_k}=v_k$, $W_{F_k}=w_k$).

In particular it is clear that
$U^2=V^2=W^2=-Id$ and that $UV=-VU=W$, $VW=-WV=U$, $WU=-UW=V$, so that
$Id,U,V$ and $W$ generate an algebra which is isomorphic to $\HH$.

\

Our space $X=X(\HH)$ is then defined inductively as the completion
of $c_{00}$ in the smallest norm satisfying the following equation:
$$\norm{x}=\norm{x}_{c_0} \vee \sup\{f(n)^{-1}\sum_{i=1}^n \norm{E_i x}:
2 \leq n, E_1<\cdots<E_n\}$$
$$ \vee \sup\{|x^*(Ex)|: k \in K, x^* \in B_k^*(X), E
\subset \N\} \vee \norm{Ux} \vee \norm{Vx} \vee \norm{Wx},$$
where $E$, and $E_1,\ldots,E_n$ are intervals of integers.

\ 

We may immediately observe that $U,V$ and $W$ extend to norm $1$
operators on $X$, and even isometries on $X$, by the quaternionic
relations between them. Note also that $U, V$ and $W$ commute with any
interval projection, and that whenever $x<y$ and $T \in \{U,V,W\}$, we
have $Tx<Ty$. It follows that $\norm{Tx}_{(N)} = \norm{x}_{(N)}$, whenever $N \geq 1$
and $T \in \{U,V,W\}$; when a sequence $x_1,\ldots,x_N$ satisfies the RIS($\epsilon$)
condition, then so does $Tx_1,\ldots,Tx_N$. The adjoints
$Id_{X^*},W^*,V^*,U^*$, in this order, also satisfy the quaternionic
commutation relations, the commutation with interval projections, and the
relation with successive vectors. It follows that when
$x^* \in A^*_m$ for some $m \in \N$ and $T \in \{U,V,W\}$, we have
$T^*x^* \in A^*_m$. However, and this is fundamental,
the sequence $T^*x_1^*,\ldots,T^*x_k^*$ is not in general special
 when
$x^*_1,\ldots,x_k^*$ is a special sequence.

The next lemma is taken directly from \cite{GM2}.

\begin{lemm}\label{gm6} For
every  $K_0 \subset K$, there is a function
 $g \in \cal F$ such that $f \geq g \geq {f}^{1/2}$, $g(k)={f(k)}^{1/2}$ whenever $k \in
K_{0}$, and $g(x)=f(x)$ whenever $N \in J \setminus K_{0}$ and $x$ in the interval $[\log N,
\exp N]$.
\end{lemm}

\begin{lemm}\label{gm7} Let $0<\epsilon \leq 1, 0 \leq \delta <1$, $M \in L$ and let $n$
and $N$ be integers such that $N/n \in [\log M,\exp M]$ and $f(N) \leq (1+\delta)f(N/n)$.
Assume that $x_1,\ldots,x_N$ satisfies the RIS($\epsilon$) condition and let
$x=x_1+\ldots+x_N$. Then $\norm{(f(N)/N)x}_{(n)} \leq (1+\delta)(1+3\epsilon)$. In
particular, if $n=1$, we have $\norm{(f(N)/N)x} \leq 1+3\epsilon$.
\end{lemm}
\pf We may reproduce the proof of Lemma 7 from \cite{GM2}, provided we
show that if a vector
$Ex$ is such that $\norm{Ex}>1$, then it is normed by a $(k,g)$-form,
where $g$ is given by Lemma \ref{gm6} in the case $K_0=K\setminus\{k\}$. To see this,
note that the only new case with respect to the classical Gowers-Maurey's
proof is if $Ex$ was normed by some $T_1^*E_1^* \ldots T_m^*E_m^*x^*$, with $T_i
\in
\{U,V,W\}$ and $E_i$ an interval projection, for all $1 \leq i \leq m$,
and
$x^*$ a
$(k,g)$-form. By the commutation properties of $U,V,W$, we may assume
$Ex$ is normed by $T^*x^*$ with $T \in \{U,V,W\}$ and $x^*$ a
$(k,g)$-form. But in this case,
$T^*x^*$ is also a $(k,g)$-form, since $T^*$ is an isometry which
respects successive vectors. \pff

\

To reproduce the proof of Gowers and Maurey, after having added the
isometries $U$, $V$ and $W$ in the definition of the norm, we shall need to distinguish
the action of a functional $x^*$ from the actions of $U^*x^*$, $V^*x^*$ and $W^*x^*$. 
This is expressed by the next lemma.

\begin{lemm}\label{hahnbanach} Let $x$ be a finitely supported
vector in
$X$. Then there exists a functional $x^*$ of norm at most $1$ such that
$x^*(x) \geq 1/2\norm{x}$ and such that $x^*(Ux)=x^*(Vx)=x^*(Wx)=0$.
\end{lemm}

\pf We observe that for any reals $\alpha,\beta,\gamma$,
the inverse of the operator $Id-\alpha U - \beta V -\gamma W$ is equal to
$(1+\alpha^2+\beta^2+\gamma^2)^{-1}(Id+\alpha U+\beta V+\gamma W)$.
It follows that $$\norm{(Id-\alpha U - \beta V -\gamma W)^{-1}} \leq
\frac{1+|\alpha|+|\beta|+|\gamma|}{1+\alpha^2+\beta^2+\gamma^2} \leq
3/2$$ by elementary calculus. So for any $x \in c_{00}$,
$$d(x,[Ux,Vx,Wx])
\geq 2/3 \norm{x}.$$ We conclude using the Hahn-Banach theorem.\pff

\

Given $N \in L$, and $\delta>0$, define a {\em $\delta$-norming $N$-pair} to be a pair
$(x,x^*)$ defined as follows. Let $y_1,\ldots,y_N$ be a sequence
satisfying the RIS(1) condition. Let
$x=N^{-1}f(N)(y_1+\ldots+y_N)$. Let, for $1
\leq i
\leq N$,
$y_i^*$ be a
 functional of norm at most $1$, such that
$ran(y_i^*) \subset ran(y_i)$ and $y_i^*(y_i)=\delta$.
Let $x^*=f(N)^{-1}(y_1^*+\ldots+y_N^*)$. 
 Note that if
$(x,x^*)$ is a $\delta$-norming
$N$-pair, then $x^*(x)=\delta$ and
Lemma \ref{gm7} implies that $\norm{x}_{(\sqrt{N})} \leq 8$.
By Lemma \ref{gm4} and Hahn-Banach Theorem, $\delta$-norming
$N$-pairs $(x,x^*)$ with arbitrary constant $\delta \leq 1/2$ exist with $x$ in an
arbitrary block-subspace and
$N$ arbitrary. 

\begin{prop} The space $X(\HH)$ is hereditarily indecomposable.
\end{prop}

\pf Write $X=X(\HH)$. Let $Y$ and $Z$ be subspaces of $X$ and
$\epsilon>0$. We may assume that $Y$ and $Z$ are generated by successive
vectors in
$X$. Let $k \in K$ be such that $(\epsilon/72)f(k)^{1/2}>1$. We construct
sequences $x_1,\ldots,x_k$ and $x^*_1,\ldots,x_k^*$ as follows. Let
$N_1=j_{2k}$ and by Lemma \ref{gm4}, let $(x_1,x_1^*)$ be a $1/3$-norming $N_1$-pair
 such that
$x_1
\in Y$, with 
 $|x_1^*(Tx_1)| < k^{-1}$ if $T=U, V$ or $W$: this is possible
by Lemma
\ref{hahnbanach} applied to each of the $N_1$ vectors forming $x_1$. Since we allow an error
term
$k^{-1}$, $x_1$ and the functional
$x_1^*$ may be perturbed so that $x_1^*$ is in
$\bf Q$ and
$\sigma(x_1^*)>f^{-1}(4)$. In general, after the first $i-1$ pairs were
constructed, let $(x_i,x_i^*)$ be a $1/3$-norming $N_i$-pair such that $x_i$ and
$x_i^*$ are supported after $x_{i-1}$ and $x_{i-1}^*$, with $x_i \in Y$
if $i$ is odd and $x_i \in Z$ if $i$ is even, such that
 $|x_i^*(Tx_i)| <k^{-1}$ if $T=U, V$ or $W$,
 having perturbed
$x_i^*$ in such a way that $N_{i+1}=\sigma(x_1^*,\ldots,x_i^*)$ satisfies
$f(N_{i+1})>2^{i+1}$ and $\sqrt{f(N_{i+1})}>2|ran(\sum_{j=1}^i x_j)|$.
Now let $y=x_1+x_3+\ldots+x_{k-1}$, $z=x_2+x_4+\ldots+x_k$. Let
also $x^*=f(k)^{-1/2}(x_1^*+\ldots+x_k^*)$. Our construction guarantees
that $x^*$ is a special functional, and therefore of norm at most $1$. It
is also clear that
$y
\in Y$,
$z
\in Z$, and that $$\norm{y+z} \geq x^*(y+z) \geq 1/3 kf(k)^{-1/2}.$$

Our aim is now to obtain an upper bound for $\norm{y-z}$. Let $x=y-z$.
Let
$g$ be the function given by Lemma \ref{gm6} in the case
$K_0=K\setminus\{k\}$. By the definition of the norm, all vectors $Ex$
are either normed by
$(M,g)$-forms, by special functionals of length $k$, by images of such
functionals by $U$, $V$ or $W$, or they have norm at most $1$. In order
to apply Lemma \ref{gm3}, it is enough to show that
$|T^*z^*(Ex)|=|z^*(TEx)| \leq 1$ for any special functional $z^*$ of
length $k$ in $K$, $E$ an interval, $T$ in the set $\{Id,U,V,W\}$. Let
$z^*=f(k)^{-1/2}(z_1^*+\ldots+z_k^*)$ be such a functional with $z_l^*
\in A^*_{m_l}$ for $1 \leq l \leq k$.

We evaluate $|z_l^*(ETx_i)|$ for $1 \leq l \leq k$ and $1 \leq i \leq k$.
Recall that $T$ and $E$ commute.

 Let $t$ be the largest integer such that $m_t=N_t$. Then
$z_i^*=x_i^*$ for all $i<t$.   There are at most two values of $i<t$ such
that
$x_i \neq Ex_i \neq 0$ or $z_i^* \neq Ez_i^* \neq 0$, and for them
$|z_i^*(ETx_i)| \leq 1$.
The values of $i<t$ for which $x_i=Ex_i$ and
$z_i^*=Ez_i^*$ form an interval $e$ and satisfy
$z_i^*(Tx_i)=x_i^*(x_i)=1/3$ if
$T=Id$, or $|z_i^*(Tx_i)|=|x_i^*(Tx_i)| \leq k^{-1}$,  when $T=U,V$ or
$W$.   Therefore $|\sum_{i \in e}
z_i^*(ET(-1)^i x_i)| \leq 1$.

When $i=l=t$, we obtain $|z_t^*(TEx_t)| \leq 1$.

If $i=l>t$ or $i \neq l$ then
$z_l^*(Tx_i)=(T^*z_l^*)(x_i)$ and we have $T^*z_l^* \in A^*_{m_i}$ for
some $m_l$. Moreover, because $\sigma$ is injective and by definition of
$t$, $m_l \neq N_i$. If $m_l<N_i$, then by the remark after the
definition of $N$-pairs,
$\norm{x_i}_{\sqrt{N_i}} \leq 8$, so the lower bound of $j_{2k}$ for
$m_1$ tells us that $|T^*z_l^*(x_i)| \leq k^{-2}$. If $m_l>N_i$ the same
conclusion follows from Lemma \ref{gm1}.
There are also at most two pairs $(i,l)$ for which $0 \neq z_l^*(ETx_i)
\neq z_l^*(Tx_i)$, in which case $|z_l^*(ETx_i)| \leq 1$.

Putting these estimates together we obtain that
$$|z^*(TEx)|=|z^*(TE(\sum_{i=1}^k (-1)^i x_i)| \leq
f(k)^{-1/2}(2+1+1+2+k^2.k^{-2}) \leq 1.$$

We also know that $(1/8)(x_1,\ldots,x_k)$ satisfies the RIS(1) condition.
Hence, by Lemma \ref{gm3}, $\norm{x} \leq 24 kg(k)^{-1}= 24kf(k)^{-1}$.
It follows that $\norm{y-z} \leq 72f(k)^{-1/2} \norm{y+z} \leq \epsilon
\norm{y+z}$. It follows that $Y$ and $Z$ do not form a direct sum and so
$X$ is H.I..
\pff

\

We may also construct a complex version $X(\C)$ of our space, with a
$2$-dimen\-sio\-nal decomposition, and a canonical isometry $J$ satisfying
$J^2=-Id$ corresponding to a representation of the complex numbers as operators
on each $2$-dimensional summand. We leave the reader adapt our
definitions and proofs to that case. Alternatively one could use the
previous $4$-dimensional decomposition setting and put only the
operator $U$, instead of $U,V$ and $W$, in the definitions and the
proofs. We therefore obtain:

\begin{prop} The space $X(\C)$ is hereditarily indecomposable.
\end{prop}

\section{Properties of operators on $X(\HH)$ and on $X(\C)$}

We now turn to the operator properties of our spaces
$X(\C)$ and $X(\HH)$. The quarternionic case is immediate from Theorem \ref{valentin}.

\begin{prop} Let $X=X(\HH)$. Then the algebra ${\cal L}(X)/{\cal
S}(X)$ is isomorphic to $\HH$. Furthermore, for any $Y \subset X$,
$\dim {\cal L}(Y,X)/{\cal S}(Y,X)=4$, i.e. every operator from $Y$ into
$X$ is of the form $a i_{Y,X}+bU_{|Y}+cV_{|Y}+dW_{|Y}+s$, where $a,b,c,d$
are reals and
$s$ is strictly singular.
\end{prop}

\pf The operators $Id$, $U$, $V$, and $W$ generate a quaternionic
division algebra, so ${\cal L}(X)/{\cal S}(X)$ is of dimension at least
$4$. By Theorem \ref{valentin}, it is isomorphic to $\HH$, and
furthermore, since
${\cal L}(X)/{\cal S}(X)$ embeds into $E_Y={\cal L}(Y,X)/{\cal S}(Y,X)$,
and $E_Y$ is of dimension at most $4$ for any $Y \subset X$ by \cite{F2},
we deduce that $\dim E_Y=4$ for any $Y \subset X$.
\pff

\

The complex case requires the following lemma, which is inspired by 
Lemma 4.14 from
\cite{A}.

\begin{lemm} \label{sss} Let $X$ be a real H.I. space and $J$ be an
operator on
$X$ such that $J^2=-Id$. Let
$Y$ be a subspace of
$X$.  Let $T \in
{\cal L}(Y,X)$ be an operator which is not of the form $\lambda
i_{Y,X}+\mu J_{|Y} +s$, with $\lambda$, $\mu$ scalars and $s$ 
strictly singular. Then there exists a finite codimensional subspace
$Z$ of $Y$ and some $\alpha>0$ such that:
$$\forall z \in Z, d(Tz,[z,Jz]) \geq \alpha \norm{z}.$$
\end{lemm}

\pf Otherwise we may construct a normalized basic sequence $(y_n)
\in Y$, and scalars sequences $(\lambda_n), (\mu_n)$ with for all $n \in
\N$,
$$\norm{Ty_n -\lambda_n y_n-\mu_n Jy_n} \leq 2^{-n}.$$
It follows that for $C=1+\norm{T}$, for all $n \in \N$,
$$\norm{\lambda_n y_n+\mu_n Jy_n} \leq C.$$
We may assume for convenience that
$J$ is isometric. We note that
$$(\lambda_n Id+\mu_n
J)^{-1}=\frac{\lambda_n Id-\mu_n J}{\lambda_n^2+\mu_n^2},$$
from which it follows that 
$$1=\norm{y_n} \leq C \frac{\norm{\lambda_n Id-\mu_n
J}}{\lambda_n^2+\mu_n^2},$$   
 so $$\lambda_n^2+\mu_n^2 \leq C(|\lambda_n|+|\mu_n|).$$
It follows immediately that $\max(|\lambda_n|,|\mu_n|) \leq 2C$. Thus
we may assume that the sequences $(\lambda_n)$
and $(\mu_n)$ converge, and, passing to a subsequence, deduce that for
some
$\lambda$,
$\mu$,
$$\norm{Ty_n-\lambda y_n-\mu Jy_n} \leq 3.2^{-n}.$$
From this it follows that the restriction of $T-\lambda i_{Y,X} - \mu
J_{|Y}$ to the space generated by the basic sequence $(y_n)$ is compact,
therefore strictly singular. By Lemma \ref{quasimaximal}, we deduce that
$T-\lambda i_{Y,X} - \mu J_{|Y}$ is strictly singular on $Y$, a
contradiction. \pff

\begin{prop}\label{operatorscomplex} Let $X=X(\C)$. Then the algebra
${\cal L}(X)/{\cal S}(X)$ is isomorphic to $\C$. Furthermore, for any $Y
\subset X$,
$\dim {\cal L}(Y,X)/{\cal S}(Y,X)=2$, i.e. any operator from $Y$ into $X$
is of the form $\lambda i_{Y,X}+\mu J_{|Y}+s$, where $\lambda,\mu$ are
reals and $s$ is strictly singular.
\end{prop}

\pf Let $Y$ be a subspace of $X$. The operator $J_{|Y}$ is not of the
form $\lambda i_{Y,X}+s$, $s$ strictly singular, otherwise by the H.I.
property, Lemma \ref{quasimaximal}, $J-\lambda Id$ would be strictly
singular, which would contradict the fact that $J^2=-Id$. It follows that
$\dim {\cal L}(Y,X)/{\cal S}(Y,X) \geq 2$.

We assume  $\dim {\cal L}(Y,X)/{\cal S}(Y,X) >2$ and look for a
contradiction. Let $T \in {\cal L}(Y,X)$ which is not of the form 
$\lambda i_{Y,X}+\mu J_{|Y}+s$ and assume without loss of generality that
$\norm{T} \leq 1$. Then by Lemma
\ref{sss} we may find some $\alpha>0$ and some 
subspace $Z$ of $Y$ such that for all $z \in Z$,
$$d(Tz,[z,Jz]) \geq \alpha\norm{z}.$$
We may assume that $Z$ is generated by successive vectors  with respect
to the
$2$-dimensional decomposition of
$X$.

We fix a sequence $(y_n)$ in $Z$ such that for all $n$, $y_{n+1}$ and
$Ty_{n+1}$ are supported after $y_n$ and $Ty_n$, $\norm{y_n}_{(n)} \leq 1$
while $\norm{y_n} \geq 1/2$. Let $k \in K$ and construct sequences
$x_1,\ldots,x_k$ and
$x_1^*,\ldots,x_k^*$ as follows. Let $N_1=j_{2k}$.
Let $x_1=N_1^{-1}f(N_1)(y_{n_1}+\ldots+y_{n_{N_1}})$ where $y_{n_i}$ is a
subsequence satisfying the RIS(1) condition. Then 
$Tx_1=N_1^{-1}f(N_1)(Ty_{n_1}+\ldots+Ty_{n_{N_1}})$, where
the sequence $(Ty_{n_i})$  satisfies the RIS(1)
condition as well. We let
$x_1^*=f(N_1)^{-1}(y_{n_1}^*+\ldots+y_{n_{N_1}})$ be associated to
$Tx_1$ so that
$(x_1^*,Tx_1)$ is an $\alpha/2$-norming $N_1$-pair, and such that 
$|x_1^*(x_1)|<k^{-1}$ and $|x_1^*(Jx_1)|<k^{-1}$ (to find $y_{n_i}^*$'s realizing this, apply
the inequality from Lemma
\ref{sss} to each $y_{n_i}$, with Hahn-Banach theorem). Lemma
\ref{gm7} implies that $\norm{Tx_1}_{(\sqrt{N_1})} \leq 8$.
Repeating this, and up to perturbations,
we build, for $1 \leq i \leq k$, $x_i \in Z$ and $x_i^*$ so that
$(x_i^*,Tx_i)$ is an
$\alpha/2$-norming $N_i$-pair, and such that $x_1^*,\ldots,x_k^*$ is
a special sequence. We let $x=x_1+\ldots+x_k$ and
$x^*=f(k)^{-1/2}(x_1^*+\ldots+x_k^*)$. We therefore have
$$\norm{x} \geq \norm{Tx} \geq x^*(Tx) \geq (\alpha/2) kf(k)^{-1/2}.$$

We now use Lemma \ref{gm3} to obtain an upper estimate for $\norm{x}$.  
Let
$g$ be the function given by Lemma \ref{gm6} in the case
$K_0=K\setminus\{k\}$. By the definition of the norm, all vectors $Ex$
are either normed by
$(M,g)$-forms, by special functionals of length $k$, by images of such
functionals by $J$, or they have norm at most $1$. In order
to apply Lemma \ref{gm3}, it is enough to show that
$|z^*(Ex)| \leq 1$
and $|J^*z^*(Ex)|=|z^*(JEx)| \leq 1$ for any special functional $z^*$ of
length $k$ in $K$ and any interval $E$. Let
$z^*=f(k)^{-1/2}(z_1^*+\ldots+z_k^*)$ be such a functional with $z_j^*
\in A^*_{m_j}$.

We evaluate $|z_l^*(Ex_i)|$ and $|z_l^*(EJx_i)|$ for $l \leq k$ and $i
\leq k$.

 Let $t$ be the largest integer such that $m_t=N_t$. Then
$z_i^*=x_i^*$ for all $i<t$.   There are at most two values of $i<t$ such
that
$x_i \neq Ex_i \neq 0$ or $z_i^* \neq Ez_i^* \neq 0$, and for them
$|z_i^*(Ex_i)| \leq 1$ and $|z_i^*(JEx_i)| \leq 1$.
The values of $i<t$ for which $x_i=Ex_i$ and
$z_i^*=Ez_i^*$ give $|z_i^*(Ex_i)|=|x_i^*(x_i)|<k^{-1}$ and 
$|z_i^*(EJx_i)|=|x_i^*(Jx_i)|<k^{-1}$.

When $i=l=t$, we obtain $|z_t^*(Ex_t)| \leq 1$ and $|z_t^*(EJx_t)| \leq
1$.

If $i=l>t$ or $i \neq l$ then
$z_l^*(Jx_i)=(J^*z_l^*)(x_i)$ and we have as before $J^*z_l^* \in
A^*_{m_l}$ for some $m_l \neq N_i$. If $m_l<N_i$, then as we remarked
above,
$\norm{x_i}_{\sqrt{N_i}} \leq 8$, so the lower bound of $j_{2k}$
for
$m_1$ tells us that $|J^*z_l^*(x_i)| \leq k^{-2}$. If $m_l>N_i$
the same conclusion follows from Lemma \ref{gm1}.
There are also at most two pairs $(i,l)$ for which $0 \neq z_l^*(EJx_i)
\neq z_l^*(Jx_i)$, in which case $|z_l^*(EJx_i)| \leq 1$.
The same proof holds for $|z_l^*(Ex_i)|$.

Putting these estimates together we obtain that
$$|z^*(Ex)| \vee |z^*(EJx)| \leq
f(k)^{-1/2}(2+k.k^{-1}+1+2+k^2.k^{-2}) \leq 1.$$

We also know that $(1/8)(x_1,\ldots,x_k)$ satisfies the RIS(1)
condition. Hence, by Lemma \ref{gm3}, $\norm{x} \leq 24
kg(k)^{-1}= 24kf(k)^{-1}$. 

\

Finally we deduce that $\alpha\sqrt{f(k)} \leq 48$, a
contradiction, since $k$ was arbitrary in $K$.
We conclude that $\dim {\cal L}(Y,X)/{\cal S}(Y,X)=2$ for any $Y \subset
X$, and that ${\cal L}(X)/{\cal S}(X)$ is isomorphic to $\C$.
\pff

\section{Complex structures on $X(\C)$ and $X(\HH)$}

Given $X$ a real Banach space, let ${\cal GL}(X)$ denote the group of invertible operators on $X$, and
 let ${\cal I}(X)$ denote the subset of operators $I$ on $X$ such that $I^2=-Id$.

\begin{lemm}\label{perturbation} Let $X$ be a real Banach space. Let $I \in {\cal I}(X)$, and let $S \in {\cal S}(X)$
be such that $I+S \in {\cal I}(X)$. Then the complex structures $X^I$ and $X^{I+S}$ associated to $I$ and $I+S$
respectively are isomorphic.
\end{lemm}

\pf Write $T=I+S$. From $T^2=-Id$, we deduce $S^2=-IS-SI$. Now let
$\alpha=Id-(SI/2)$. This is an $\R$-linear map on $X$.
Moreover, it is easy to check, using the relation satisfied by $S^2$,
that
$$\alpha I=I +(S/2)=(I+S)\alpha=T\alpha.$$
This relation  ensures that $\alpha$ may be seen as a
$\C$-linear operator from $X^{I}$ into $X^{T}$.

Furthermore, by properties of strictly singular operators (see e.g.
\cite{LT}), the $\R$-linear operator $\alpha$  is
Fredholm of index
$0$ on
$X$ as a strictly singular perturbation of $Id$.
This means that $\alpha(X)$ is closed, and that 
$$\dim_{\R} Ker \alpha = \dim_{\R} (X/\alpha(X)) <+\infty.$$
This is also true when $\alpha$ is seen as $\C$-linear (note in
particular that $Ker \alpha$ is $I$-stable and $\alpha(X)$ is
$T$-stable). That is, $\alpha(X^I)$ is  closed in $X^T$,
and 
$$\dim_{\C} Ker \alpha = \dim_{\C} (X^T/\alpha(X^I)) <+\infty.$$

 Therefore $\alpha$ is $\C$-Fredholm with index $0$, i.e. there exist
$\C$-linear decompositions
$X^I=X_0 \oplus F_0$ and
$X^{T}=Y_0 \oplus G_0$, with $\dim_{\C} F_0=\dim_{\C} G_0<+\infty$,
such
that the restriction of $\alpha$ to $X_0$ is a $\C$-linear isomorphism
onto
$Y_0$.
Since $F_0$ and $G_0$ are isomorphic, we deduce that there exists a
$\C$-linear isomorphism from $X^I$ onto $X^T$. \pff

\
 
 Let $\pi$ denote the quotient map
from
${\cal L}(X)$ onto ${\cal L}(X)/{\cal S}(X)$.  We also let
$({\cal L}(X)/{\cal S}(X))_0$ denote the group $\pi({\cal GL}(X))$, and
 $\tilde{{\cal I}}(X)$ denote the set of elements of $({\cal L}(X)/{\cal S}(X))_0$ whose square is
equal to $-\pi(Id)$. In the following we shall identify a complex structure on $X$
with the associated operator $I
\in {\cal I}(X)$.

\begin{prop}\label{principe} Let $X$ be a real Banach space. Then the quotient map $\pi$ induces an injective map
$\tilde{\pi}$ from the set of isomorphism classes of complex structures on $X$ into the set of conjugation classes of
elements of
$\tilde{{\cal I}}(X)$ for the group $({\cal L}(X)/{\cal S}(X))_0$. The image of $\tilde{\pi}$ is the set
of conjugation classes of elements of $\tilde{{\cal I}}(X)$ which may be lifted to an element of ${\cal I}(X)$.
If ${\cal S}(X)$ admits a supplement in ${\cal L}(X)$ which is a subalgebra of ${\cal L}(X)$,
then $\tilde{\pi}$ is bijective.
\end{prop}

\pf For any operator $T$ on $X$, we write $\tilde{T}=\pi(T)$. Let $I$ and $T$ be operators in ${\cal I}(X)$. If
$\alpha$ is a $\C$-linear isomorphism from
$X^I$ onto $X^T$, then the $\C$-linearity means that $\alpha I=T \alpha$.
Therefore $\tilde{\alpha}\tilde{I}=\tilde{T}\tilde{\alpha}$, and $\tilde{I}$ and $\tilde{T}$ satisfy a conjugation
relation. Conversely, if $\tilde{I}=\tilde{\alpha}^{-1} \tilde{T} \tilde{\alpha}$ for some $\alpha \in {\cal GL}(X)$,
then  $\alpha^{-1}T\alpha=I+S$, where $S$ is strictly singular.
Note that
 $(I+S)^2=-Id$, and since $T\alpha=\alpha(I+S)$, $\alpha$ is a $\C$-linear isomorphism from 
$X^{I+S}$ onto $X^T$.
By Lemma \ref{perturbation}, it follows that $X^I$ and $X^T$ are isomorphic.
This proves that $\tilde{\pi}$ is well-defined and injective.

If ${\cal H}(X)$ is a subalgebra of ${\cal L}(X)$ which supplements ${\cal S}(X)$, then let $T \in {\cal L}(X)$ be
such that $\tilde{T}^2=-\tilde{Id}$; we may assume that $T$ (and therefore $T^2$) belongs to ${\cal H}(X)$. Then
since
$T^2=-Id+S$,
$S$ strictly singular, $T^2$ must be equal to $-Id$. Any class $\tilde{T} \in \tilde{\cal I}(X)$ may  therefore be
lifted to an element of ${\cal I}(X)$. \pff

\

We recall that $X^J(\C)$ denotes the complex structure on $X(\C)$ associated to 
the operator $J$, and that
two Banach spaces $X$ and $X'$
are said to be totally incomparable if no subspace of $X$ is isomorphic
to a subspace of
$X'$.

\begin{prop} \label{XJ} The space $X^J(\C)$ and its conjugate $X^{-J}(\C)$ are
  complex H.I. and totally incomparable. Moreover, any complex
structure on $X(\C)$ is isomorphic either to $X^J(\C)$ or to $X^{-J}(\C)$.
\end{prop}
\pf Any complex structure on $X(\C)$ is H.I., since $X(\C)$ is real H.I..
We have ${\cal L}(X(\C))=[Id,J] \oplus {\cal S}(X(\C))$, and 
$$({\cal L}(X(\C))/{\cal S}(X(\C)))_0 \simeq \C^*.$$
The only two elements of $\C$ of square $-1$ are $i$ and $-i$. By Proposition~\ref{principe}, it follows that $X^J(\C)$
and
$X^{-J}(\C)$ are the only two complex structures on $X(\C)$ up to isomorphism.

Assume now $\alpha$ is a $\C$-linear map from a $\C$-linear
subspace $Y$ of
$X^J(\C)$ into
$X^{-J}(\C)$. This is in particular an $\R$-linear map from $Y$ into
$X(\C)$.  So by Proposition
\ref{operatorscomplex},
$\alpha=a Id_{|Y} +bJ_{|Y}+s$,
where
$s$ is strictly singular. The fact
that $\alpha$ is
$\C$-linear means that
$\alpha J_{|Y}=-J\alpha$. This implies by ideal properties of strictly singular operators that
 $a Id_{|Y}+bJ_{|Y}$ is strictly singular and therefore, $a=b=0$. It follows
that $\alpha$ is
$\R$-strictly singular and thus $\C$-strictly singular. Therefore $X^{J}(\C)$
and
$X^{-J}(\C)$ are totally incomparable.
\pff

\

Following \cite{GM1}, we shall say that an operator $W \in {\cal L}(Y,Z)$
is {\em finitely singular} if the restriction of $W$ to some
finite-codimensional subspace of $Y$ is an isomorphism into $Z$. This
means that $WY$ is closed and that the Fredholm index $i(W)$ is defined
with values in
$\Z \cup \{-\infty\}$, where
$$i(W)=\dim(Ker(W))-\dim(Z/WY).$$ 
The next proposition and corollaries show that the example of $X(\C)$ is essentially the only
one to ensure the total incomparability property.

\

\begin{prop}\label{strange}
Let $X$ be a real Banach space, $T,U \in {\cal I}(X)$. Then $X^T$ is isomorphic to
$X^U$, or $T-U$ induces a $\C$-linear isomorphism
from a subspace of $X^T$ onto a subspace of $X^{-U}$. 
\end{prop}

\pf We note that $(T+U)T=U(T+U)$, therefore $T+U$ is $\C$-linear from
$X^T$ into $X^U$. The similar result
holds for $T-U$ between $X^T$ and $X^{-U}$.

Assume $T-U$ does not induce a $\C$-linear isomorphism
from a subspace of $X^T$ onto a subspace of $X^{-U}$. Then $T-U$ is strictly singular
on $X^T$, and we intend to deduce that $T+U$ is "essentially" an isomorphism from
$X^T$ onto $X^U$.

\

First we note that $T+U$ is finitely singular on $X$. To see this it is enough
by de\-fi\-nition to prove that $T+U$
is finitely singular as a $\C$-linear operator from $X^T$ into
$X^U$. If this were false, then by \cite{LT} Proposition 2.c.4, we could
find a ($\C$-linear) infinite dimensional subspace $Y$ of $X^T$ such that
$\norm{(T+U)_{|Y}} < \norm{T}^{-1}$. Since $T-U$ is strictly singular on $X^T$, 
we would find a norm $1$ vector $y$ in $Y$ with $\norm{(T-U)y} < \norm{T}^{-1}$. We 
would then deduce that
$$\norm{Ty} \leq 2^{-1}(\norm{(T+U)y}+\norm{(T-U)y}) <\norm{T}^{-1}\norm{y},$$
a contradiction.

\

We now prove that whenever $\lambda \in ]0,1[$, the operator
$T+\lambda U \in {\cal L}(X)$ is finitely singular. Assume on the contrary that
$T+\lambda U$ is not finitely singular for some $\lambda \in ]0,1[$. Let
$c=(1-\lambda^2)(2(1+\norm{T}+2\lambda\norm{U}))^{-1}$ and take an
arbitrary
$0<\epsilon<c$. As before there exists some infinite-dimensional ($\R$-linear)
subspace $Y$ of $X$ such that $\norm{(T+\lambda U)_{|Y}}<\epsilon$.
Therefore, for all $y \in Y$,
$$ \norm{Ty+\lambda Uy} \leq \epsilon\norm{y} \leq c\norm{y}, \leqno (1)$$
and by composing by $U$,
$$ \norm{UTy-\lambda y} \leq \epsilon\norm{U}\norm{y} \leq
c\norm{U}\norm{y}. \leqno (2)$$ We deduce from this that $Y$ and $TY$ form a direct
sum. Indeed, if there existed norm $1$ vectors $y$ and $z$ in $Y$ with
$\norm{z-Ty} \leq c$, then we would have, by (1), $\norm{z+\lambda Uy}
\leq 2c$, therefore
$$ \norm{Uz-\lambda y} \leq 2c\norm{U}, \leqno (3)$$
but also $\norm{Tz+y} \leq c\norm{T}$, so by (1) again,
$$ \norm{y-\lambda Uz} \leq c(1+\norm{T}). \leqno (4)$$
From (3) and (4), we would get
$$1-\lambda^2=\norm{(1-\lambda^2)y} \leq c(1+\norm{T}+2\lambda\norm{U}),$$
a contradiction by choice of $c$. Therefore $Y+TY$ forms a direct sum with
projection constants depending only on $\lambda$, $\norm{T}$ and 
$\norm{U}$. 

Now since $T-U$ is $\C$-strictly singular from $X^T$ into $X^{-U}$,
and $Y \oplus TY$ is a $\C$-linear subspace of $X^T$, there exist
$y,z \in Y$ with $\norm{y+Tz}=1$ and $\norm{(T-U)(y+Tz)} \leq \epsilon$,
which means $\norm{-z-UTz+Ty-Uy} \leq \epsilon$.
By (1) and (2), and the fact that $Y \oplus TY$ is direct, we deduce
$$\norm{-z-\lambda z+Ty+\lambda^{-1} Ty} \leq
\epsilon+\epsilon\norm{U}\norm{z}+\epsilon\lambda^{-1}\norm{y} \leq
C\epsilon,$$
where $C$ depends only on $\lambda$, $\norm{T}$ and $\norm{U}$.
In the same way,
$$\norm{-z-\lambda z +Ty+\lambda^{-1}Ty} \geq C'((1+\lambda)\norm{z}
+(1+\lambda^{-1})\norm{Ty}) \geq C'',$$
where again $C''$ depends only on $\lambda$, $\norm{T}$ and $\norm{U}$.
As $\epsilon$ was arbitrary, we obtain a contradiction.

\

We have therefore proved that $T+\lambda U$ is finitely singular
whenever $\lambda \in ]0,1]$, and this is obvious for $\lambda=0$. In
other words, the Fredholm index of $T+\lambda U$ is defined
for all $\lambda \in [0,1]$. By continuity of the Fredholm index, 
\cite{LT} Proposition 2.c.9, we deduce that $ind(T+U)=ind(T)=0$,
since $T$ is an isomorphism. Therefore $T+U$ is Fredholm with index $0$.
It is therefore also Fredholm with index $0$ as a $\C$-linear operator
from $X^T$ into $X^U$, and we deduce that $X^T$ and $X^U$ are isomorphic.
\pff

\begin{coro}\label{ttt}
Let $X$ be a real Banach space with two totally incomparable complex structures. Then
these complex structures are conjugate up to isomorphism and
both saturated with H.I. subspaces.
\end{coro}

\pf Assume $X^T$ is totally incomparable with $X^U$. By Proposition \ref{strange} (applied to
$U$ and
$-T$), $X^U$ is isomorphic to $X^{-T}$.
To show that $X^T$ is HI-saturated, it is enough by Gowers' dichotomy theorem to prove
that $X^T$ does not contain a subspace with an unconditional basis. Indeed,
if $Y$ were such a subspace, then by the remark in the introduction, $Y$ 
would be isomorphic to $\overline{Y}$, $\C$-linear subspace of $X^{-T}$, 
 which would contradict the total incomparability of $X^T$ with
$X^{-T}$. \pff

\begin{coro}
There cannot exist more than two mutually totally incomparable complex structures on
a Banach space.
\end{coro}

Note that  $X^{J}(\C) \oplus X^{J}(\C)$ is  a non H.I. space which is totally
incomparable with its conjugate
$X^{-J}(\C) \oplus X^{-J}(\C)$. Therefore there is no obvious direction in which to improve
Corollary \ref{ttt}. 

G. Godefroy asked the author whether there also existed a Banach space with exactly three
complex structures. The answers turns out to be yes.

\begin{prop}
For any $n \in \N$, the space $X(\C)^{n}$ has exactly $n+1$ complex structures up to
isomorphism.
\end{prop}

\pf We have that ${\cal L}(X(\C)^n) \simeq {\cal M}_{n}(\C) \oplus {\cal
S}(X(\C)^n)$, and 
$$({\cal L}(X(/C)^n)/{\cal S}(X(\C)^n)_0 \simeq {\cal GL}_n(\C).$$
Any complex $(n,n)$-matrix whose square is equal to $-Id_{\C^n}$ has minimal polynomial
$X^2+1$, and is therefore similar to a diagonal matrix with $i$ or $-i$'s down the diagonal; and
there are 
$n+1$ similarity classes of such matrices according to the number of $i$'s. Therefore by
Proposition
\ref{principe}, there are $n+1$ complex structures on $X(\C)^n$, up to isomorphism.
If we let $X=X^J(\C)$, these structures are isomorphic to the spaces
$X^k \oplus \overline{X}^{n-k}$, $0 \leq k \leq n$. \pff

\

We now turn to the complex structures on $X(\HH)$.
 Note that since $U^2=-Id$, the real space $X(\HH)$ may be
equipped with a complex structure associated to $U$. Let
$X=X^{U}(\HH)$ denote this complex space. As it is a real H.I.
space, it is H.I. as a complex space.

Any $\R$-linear operator $T$ on $X$ is of the form $aId+bU+cV+dW+S$, where
$S$ is strictly singular. Saying that $T$ is $\C$-linear means that
$T$ commutes with $U$, which implies that $c=d=0$. Therefore we deduce
that any $\C$-linear operator on $X$ is of the form
$aId+bU+S=(a+ib).Id+S$, as expected since $X$ is H.I.! (here we used that
any $\R$-strictly singular operator on $X$ is $\C$-strictly singular).

\

S. Szarek asked whether there exists a Banach space not isomorphic to a Hilbert space, with
unique complex structure (\cite{S2}, Pb 7.2). The following proposition answers this
question by the positive.

\begin{prop}
The space $X(\HH)$ admits a unique complex structure up to isomorphism.
\end{prop}

\pf We may write ${\cal L}(X(\HH))=[Id,U,V,W] \oplus {\cal S}(X(\HH))$, and 
$$({\cal L}(X(\HH))/{\cal S}(X(\HH)))_0 \simeq \HH^*.$$

Write the generators of $\HH$ as $\{1,i,j,k\}$.
The elements of $\HH$ of square $-1$ are of the form $r=bi+cj+dk$ with $b^2+c^2+d^2=1$. 
Any element $r$ of this form is in the conjugation class of $i$, since
$i(i+r)=(i+r)r$, for $r \neq -i$, and $ij=-ji$, for $r=-i$. Therefore by Proposition \ref{principe} all 
complex structures on $X(\HH)$ are isomorphic.
 \pff

\section{Some other spaces with unique complex structure}

We show how to construct various other examples of Banach spaces with a unique complex structure. Note that
non H.I.
 examples of real spaces on which operators are of the form $\lambda Id+S$, $S$ strictly singular, may be found e.g. in
 \cite{AM}, \cite{F3}.

\begin{prop}\label{XplusX} Let $\{X_i, 1 \leq i \leq N\}$ be a family of pairwise totally incomparable
  real Banach spaces with the $\lambda Id+S$-property. For $1 \leq i \leq N$, let $n_i \in \N$.
Then $\sum_{1 \leq i \leq N} \oplus X_i^{2n_i}$ has a unique complex structure up to isomorphism.
\end{prop}

\pf If $X$ has the $\lambda Id+S$-property, and $n \in \N$, let ${\cal M}_{2n}(Id_X)$ be the space of 
$(2n,2n)$-matrix operators on $X^{2n}$
with homothetic coefficients, which we identify with the space ${\cal M}_{2n}$ of real $(2n,2n)$-matrices. 
Let ${\cal GL}_{2n}$ denote the group of invertible real $(2n,2n)$-matrices.
We have that ${\cal L}(X^{2n})={\cal M}_{2n}(Id_X) \oplus {\cal S}(X^{2n})$ and
 $$({\cal L}(X^{2n})/{\cal S}(X^{2n}))_0 \simeq {\cal GL}_{2n}.$$
 Now any real $(2n,2n)$-matrix whose square is
$-Id_{\R^{2n}}$ is diagonalizable with
 $\C$-eigenvalues $i$ and $-i$, each with multiplicity $n$.
So any two such matrices are
$\C$-similar and therefore
$\R$-similar. By Proposition \ref{principe}, it follows that all complex structures on
$X^{2n}$ are isomorphic.

If $X$ is a direct sum $\sum_{1 \leq i \leq N} \oplus X_i^{2n_i}$, where the $X_i$'s are pairwise
 totally
incomparable, then as ${\cal L}(X_i,X_j)={\cal S}(X_i,X_j)$ whenever $i \neq j$, we have
$${\cal L}(X) \simeq (\sum_{1 \leq i \leq N}\oplus {\cal M}_{2n_i}(Id_{X_i})) \oplus
{\cal S}(X),$$ and
$$({\cal L}(X)/{\cal S}(X))_0 \simeq \Pi_{1 \leq i \leq N}  {\cal GL}_{2n_i}.$$
It follows immediately from the case $N=1$ that there is a unique conjugacy class of elements of square $-1$ for the
group $\Pi_{1 \leq i \leq N}  {\cal GL}_{2n_i}$. So all complex structures on $X$ are isomorphic.
\pff

\

All the examples considered so far fail to have an unconditional basis. Indeed for each of them the
 quotient algebra
${\cal L}(X)/{\cal S}(X)$ is finite dimensional. We now show how to construct a real Banach space $X(\DD)$ with an
unconditional basis, not isomorphic to
$l_2$, and with unique complex structure (here $\DD$ stands for "$2$-block diagonal").
This is exactly the unconditional version of
$X(\C)$, precisely in the same way as the  space of Gowers \cite{G0}, on which every
operator is the sum of a diagonal and a strictly singular operator, is the
unconditional version of Gowers-Maurey's space
$X_{GM}$. In other words, it is a space with "as few" operators as possible to ensure
the existence of an unconditional basis and of a complex structure.
It is not difficult to show that the unconditionality of a $2$-dimensional
decomposition, and the existence of a map $J$ such that $J^2=-Id$ defined on each
$2$-dimensional summand, already imply that any $2$-block diagonal operator associated
to  bounded $(2,2)$-matrices must be bounded. This will motivate the following
definition of
$X(\DD)$.
We thank B. Maurey for a  discussion which clarified this example.

\

The basis $(e_i)$ is as before the natural basis of $c_{00}(\R)$. For $k \in \N$, let
$F_k=[e_{2k-1},e_{2k}]$. Notions of successivity are taken with respect to the $2$-dimensional decomposition
associated to the $F_k$'s. Let
${\cal D}_2(c_{00})$ denote the space of
$2$-block diagonal operators on
$c_{00}$, i.e. the space of operators $T$ on $c_{00}$ such that $T(F_k) \subset F_k$ for all $k \in \N$. Any operator
in ${\cal D}_2(c_{00})$
 corresponds to a  sequence $(M_n) \in {\cal M}_2^{\N}$ of real
$(2,2)$-matrices, and will be denoted
$D(M_n)$. For  
$M \in {\cal M}_2$, we shall consider the norm $\norm{M}$, when $M$ is seen as an operator
on $l_{\infty}^2$, or sometimes $\norm{M}_2$ (resp. $\norm{M}_1$), the
euclidean norm (resp. the $l_1$-norm) on ${\cal M}_2$ identified with
$\R^4$.

\

The space $X(\DD)$ is defined inductively as the completion
of $c_{00}$ in the smallest norm satisfying the following equation:
$$\norm{x}=\norm{x}_{c_0} \vee \sup\{f(n)^{-1}\sum_{i=1}^n \norm{E_i x}:
2 \leq n, E_1<\cdots<E_n\}$$
$$ \vee \sup\{|x^*(Ex)|: k \in K, x^* \in B_k^*(X), E
\subset \N\}$$
$$\vee \sup\{\norm{D(M_n)x}: \forall n \in \N, \norm{M_n} \leq 1\},$$

where $E$, and $E_1,\ldots,E_n$ are intervals of integers, and $(M_n)$ is a sequence
of $(2,2)$-matrices.

\

From the definition we observe immediately that any $2$-block diagonal operator $D(M_n)$, where 
the sequence $(M_n)$ is bounded, extends to a bounded operator on $X(\DD)$. The space of such operators
will be denoted ${\cal D}_2(X(\DD))$. Furthermore, the
norm on each $F_k$ is the $l_{\infty}$-norm; and whenever $n \in \N$,
 and $y_k,z_k$
belong to $F_k$, with $\norm{z_k} \leq \norm{y_k}$ for all $1 \leq k \leq n$, it follows
that
$$\norm{\sum_{k=1}^n z_k} \leq \norm{\sum_{k=1}^n y_k}.$$
This is a strong unconditionality property of $X(\DD)$ from which we deduce immediately that $(e_i)_{i \in \N}$ is an
unconditional basis for
$X(\DD)$.

\begin{prop}
Any operator on $X(\DD)$ may be written $D+S$, where $D \in {\cal D}_2(X(\DD))$ is 
$2$-block diagonal and $S \in {\cal S}(X(\DD))$ is strictly singular.
\end{prop}

\pf We sketch how to reproduce a proof from \cite{GM2}. Let $T$ be an operator on
$X(\DD)$ with
$0$'s down the $2$-block diagonal. First we show that if
$(x_n)$ is a sequence of successive vectors such that $\norm{x_n}_{(n)} \leq 1$, $A_n=supp(x_n)$ and for each $n$, $B_n
\cup C_n$ is a partition of
$A_n$ in two subsets, then $C_nTB_n x_n$ converges to $0$ (see \cite{GM2} Lemma 27). The proof is based on a
construction of a special sequence in the "usual" way, similar to our proofs for $X(\C)$ and $X(\HH)$. Arbitrary
choices of signs $-1$ or $1$ in the proof of \cite{GM2} correspond to arbitrary choices of norm $1$ $(2,2)$-matrices
in our case. If
$D=D(M_n) \in {\cal D}_2(X)$, with $\norm{M_n} \leq 1$ for all $n$, then $D$ preserves
successive vectors and $\norm{D}
\leq 1$. Therefore the estimates of the end of \cite{GM2} Lemma 27 based on properties of R.I.S. vectors and
$(M,g)$-forms are still valid. The other argument based on disjointness of supports of 
$y_n=B_n x_n$ and $z_n=C_n T B_n x_n$ is immediately seen to be
preserved as well.
Corollary 28 from \cite{GM2} may then be reproduced
to obtain that $(Tx_n)$ converges to $0$.  

Finally if $T$ is an operator on $X$, and $diag(T)$ is its $2$-block diagonal part, then  
$((T-diag(T))x_n)$ converges to $0$ whenever $(x_n)$ is a successive sequence such that
$\norm{x_n}_{(n)} \leq 1$. By Lemma \ref{gm4} this implies that $T-diag(T)$ is strictly singular. \pff

\

We recall that $({\cal L}(X(\DD))/{\cal S}(X(\DD)))_0$ is defined as the group of
elements  of ${\cal L}(X(\DD))/{\cal S}(X(\DD))$ which may be lifted to an invertible
operator. Likewise,
$(l_{\infty}({\cal M}_2)/c_0({\cal M}_2))_0$ is the group of elements of
$l_{\infty}({\cal M}_2)/c_0({\cal M}_2)$ which may be lifted to an invertible element of
$l_{\infty}({\cal M}_2)$, that is, to a sequence $(M_n) \in l_{\infty}({\cal M}_2)$ of real
$(2,2)$ matrices such that $M_n$ is invertible for each $n$ and the sequence $(M_n^{-1})$ is bounded.

\

\begin{lemm}\label{groups} The algebras ${\cal L}(X(\DD))/{\cal S}(X(\DD))$ and $l_{\infty}({\cal M}_2)/c_0({\cal M}_2)$ are
isomorphic, and the groups
 $({\cal L}(X(\DD))/{\cal S}(X(\DD)))_0$ and $(l_{\infty}({\cal M}_2)/c_0({\cal M}_2))_0$ are isomorphic. 
\end{lemm}

\pf Write $X=X(\DD)$. We note that ${\cal D}_2(X) \cap {\cal S}(X)$ is equal
to the set $\{D(M_n):
\lim_{n
\rightarrow +\infty} M_n=0\}$. Indeed if $(M_n)$ converges to $0$, fix $\epsilon>0$ and
$N$ such that $\norm{M_n} \leq \epsilon$ for all
$n \geq N$. Let $Y=[F_n]_{n \geq N}$. Any $y \in Y$ may be written $y=\sum_{n
\geq N}y_n$, $y_n \in F_n$, and therefore
$$\norm{D(M_n)(y)}=\norm{\sum_{n \geq N} M_n y_n} \leq
 \epsilon  \norm{\sum_{n \geq N} y_n}=\epsilon\norm{y},$$
by the  strong unconditionality properties of the basis.
This implies that $D(M_n)$ is compact and therefore strictly singular.
Conversely, if $\norm{M_n}$ does not converge to $0$ then for some $\alpha>0$ and some
 infinite set $N$,
 $\norm{M_n} \geq \alpha$ for any $n \in N$, and let $x_n \in F_n$ be a norm $1$ vector
such that $\norm{M_n x_n} \geq
\alpha$. Let $y_n=M_n x_n$. The map $C$ on $[y_n, n \in N]$ defined by $Cy_n=x_n$ is
bounded by the strong unconditionality properties
of $X(\DD)$. Therefore the restriction of $D(M_n)$ to $[x_n, n \in N]$ is an isomorphism
with inverse $C$ and $D(M_n)$ is not strictly singular. 

We deduce from this that
$${\cal L}(X)/{\cal S}(X) \simeq {\cal D}_2(X)/({\cal D}_2(X) \cap {\cal S}(X)) \simeq
l_{\infty}({\cal M}_2)/c_0({\cal M}_2).$$

Now if $T \in {\cal L}(X)$ is invertible with inverse $T'$, write $T=D+S$ and $T'=D'+S'$, with
$D=D(M_n)$, $D'=D(N_n)$ and $S, S'$ strictly singular. From
$TT'=Id_X$ we deduce that $DD'=Id_X+s$ where $s$ is strictly singular.
Furthermore $s=DD'-Id_X$ is $2$-block diagonal, and therefore of the form $D(s_n)$ where
$s_n$ converges to $0$. Therefore from $M_n N_n=Id_{\R^2}+s_n$, we deduce that for $n$
large enough,
$M_n$ is invertible and $M_n^{-1}=N_n(Id+s_n)^{-1}$ is bounded above. Modifying the first terms of
the sequences $(M_n)$ and $(N_n)$ (up to modifying $S$ and $S'$), we may assume that this is true for all $n \in \N$.

Conversely if  $M_n$ is invertible for all $n \in \N$, and $(M_n^{-1})$ is bounded, then 
$D(M_n)$ is an invertible operator with inverse $D(M_n^{-1})$.

It follows that the elements of $({\cal L}(X)/{\cal S}(X))_0$ are those that may be lifted to
a diagonal operator $D(M_n)$ where for all $n$,  $M_n$ is invertible, and $(M_n^{-1})$ is bounded; such an operator
$D(M_n)$ corresponds canonically to an invertible element of $l_{\infty}({\cal M}_2)$, and it follows that
$$({\cal L}(X)/{\cal S}(X))_0 \simeq (l_{\infty}({\cal M}_2)/c_0({\cal
M}_2))_0.$$ \pff

\

\begin{lemm}\label{existencedunrelevement} Let $A \in {\cal M}_2$.
There exists a map $f_A$ on ${\cal M}_2$ such that whenever $A^2=-Id+r$ with $\norm{r}<1$, it follows that
$(A+f_A(r))^2=-Id$, and such that if $\norm{r}<1$ then $\norm{f_A(r)} \leq \norm{A}((1-\norm{r})^{-1/2}-1)$.
\end{lemm}

\pf For $\norm{r} \leq 1$, let $f_A(r)=A(Id-r)^{-1/2}-A$, where $(Id-r)^{-1/2}$ is
defined as an infinite series
 $Id+\sum_{n \geq 1}c_n r^n$. Note that since $r=A^2+Id$, $(Id-r)^{-1/2}$ commutes with $A$.
It follows by an elementary computation that $(A+f_A(r))^2=-Id$.
Furthermore,
$$\norm{f_A(r)} \leq \norm{A}\sum_{n \geq 1}c_n \norm{r}^n=\norm{A}((1-\norm{r})^{-1/2}-1).$$
\pff

\begin{lemm}\label{elementary} Let $M \in {\cal M}_2$ be such that $M^2=-Id$. Then there exists
$P \in {\cal GL}_2$ such that $\norm{P}_2=\norm{P^{-1}}_2 \leq \sqrt{\norm{M}_1}$ and
$P\begin{pmatrix} 0 & -1 \\ 1 & 0 \end{pmatrix}P^{-1}=M$.
\end{lemm}

\pf If $M^2=-Id$ then $M$ is of the form $\begin{pmatrix} a & b \\ c & -a \end{pmatrix}$ with
$a^2=-1-bc$. If $c>0$ put $P=c^{-1/2}\begin{pmatrix} 1 & a \\ 0 & c \end{pmatrix}$, then
$P^{-1}=c^{-1/2}\begin{pmatrix} c & -a \\ 0 & 1 \end{pmatrix}$, and
$P\begin{pmatrix} 0 & -1 \\ 1 & 0 \end{pmatrix}P^{-1}=M$.
Furthermore, $\norm{P^{-1}}_2^2=\norm{P}_2^2=c^{-1}(1+a^2+c^2)=c-b \leq \norm{M}_1$.
If $c \leq 0$ then $b>0$ and a similar proof holds.
\pff

\

\begin{prop} The space $X(\DD)$ has unique complex structure up to isomorphism.
\end{prop}
Let $G$ be the group $(l_{\infty}({\cal M}_2)/c_0({\cal M}_2))_0$ and let
$I=\{g \in G: g^2=-1\}$.By Lemma \ref{groups} and Proposition \ref{principe},
it  is enough to prove that all elements of $I$ are $G$-conjugate.

Let $g \in I$, so $g$ is the class of some $(M_n)$ which is invertible 
in $l_{\infty}({\cal M}_2)$, that
is
$\norm{M_n}_1$ and $\norm{M_n^{-1}}_1$ are bounded by some constant $C$.
Since $g \in I$, the sequence $r_n=M_n^2+Id$ converges to $0$. Let $N \in \N$ be
such that for all $n \geq N$, $\norm{r_n}<1$.
Use  Lemma \ref{existencedunrelevement} to define $N_n=M_n+f_{M_n}(r_n)$ for $n \geq N$; we have therefore
that $N_n^2=-Id$ for all $n \geq N$ and that $M_n-N_n$ converges to $0$.
For $n<N$ we just put $N_n=\begin{pmatrix} 0 & -1 \\ 1 & 0 \end{pmatrix}$.
It follows that $g$ is also the class of $(N_n)$ modulo $c_0({\cal M}_2)$.

By Lemma \ref{elementary}, there exists $P_n \in {\cal GL}_2$ such that
$\norm{P_n}_2=\norm{P_n^{-1}}_2 \leq \sqrt{C}$, and such that
$$N_n=P_n\begin{pmatrix} 0 & -1 \\ 1 & 0 \end{pmatrix}P_n^{-1}.$$
Therefore $(P_n)_{n \in \N}$ and $(P_n^{-1})_{n \in \N}$ define inverse elements $p$ and
$p^{-1}$ in
$G$. Let $j$ be the element of $I$ associated to the sequence $(j_n) \in l_{\infty}({\cal
M}_2)$ of constant value
$\begin{pmatrix} 0 & -1 \\ 1 & 0 \end{pmatrix}$. We deduce that
$$g= p j p^{-1},$$
and therefore there is a unique $G$-conjugacy class of elements of $I$.\pff

\section{Final remarks and questions}

The fact that any complex Banach space which is real H.I. is also complex
H.I. raises the following question. Does there exist a complex H.I. space
which is not H.I. when seen as a real Banach space? The answer turns out
to be yes.

\

To prove this, we consider the canonical complexification $X_{GM}
\oplus_{\C} X_{GM}$ of the real version of
Gowers-Maurey's space, i.e.
$X_{GM}
\oplus X_{GM}$ with the complex structure associated to the operator $I$ defined by
 $I(x,y)=(-y,x)$. Note that by Proposition \ref{XplusX}, any other complex structure on 
$X_{GM} \oplus X_{GM}$ would be iso\-mor\-phic.

\begin{prop} The complexification of  $X_{GM}$ is complex
H.I. but not real H.I..
\end{prop}

\pf Seen as a real space, $X=X_{GM} \oplus_{\C} X_{GM}$ is clearly not
H.I. as a direct sum of infinite dimensional spaces. Let now $Y$ be a
$\C$-linear subspace of
$X$. We denote by $p_1$ and $p_2$ the projections on the first and the
second summand of $X=X_{GM} \oplus X_{GM}$ respectively. Either
$p_{1|Y}$ or $p_{2|Y}$ is not strictly singular, and without loss of
generality this is true of $p_{1|Y}$; so $p_{1|Y_1}$ is an isomorphism for
some $\R$-linear subspace $Y_1$ of $Y$. We may therefore find a
subspace
$Z$   of $X_{GM}$ and a map $\alpha: Z \rightarrow X_{GM}$ such that
$Y_1=\{(z,\alpha z), z \in Z\} \subset Y$ (take $Z=p_1(Y_1)$ and
$\alpha=p_2(p_{1|Y_1})^{-1}$). By
$\C$-linearity,
$Y_2:=iY_1=\{(-\alpha z,z), z \in Z\}$ is also an $\R$-linear subspace of
$Y$.

By the properties of $X_{GM}$, $\alpha$ is of the form
$\lambda i_{Z,X_{GM}}+s$, where $s$ is strictly singular.
For our computation we may assume that the norm $|||.|||$ on $X \oplus X$ is the $l_1$-sum norm.
Let $\epsilon>0$ be such that
$2(1+|\lambda|+\epsilon)(1+|\lambda|)\epsilon<1$. Passing to a subspace,
we may assume that $\norm{s} \leq \epsilon$.
We prove that whenever $y_1 \in Y_1$, $y_2 \in Y_2$ are norm $1$ vectors,
we have $|||y_1-y_2|||>\epsilon$. Indeed, otherwise 
let $y_1=(z_1,\alpha z_1)$ be of norm $1$ with $z_1 \in Z$, and
$y_2=(-\alpha z_2,z_2)$ be of norm $1$ with $z_2 \in Z$, with
$$\epsilon \geq |||(z_1+\alpha z_2,\alpha z_1-z_2)|||,$$
therefore
$$\epsilon \geq \norm{z_1+\lambda z_2 +s z_2},$$
and since $\norm{s} \leq \epsilon$,
$$\norm{z_1+\lambda z_2} \leq 2\epsilon.$$
Likewise,
$$\norm{\lambda z_1-z_2} \leq 2\epsilon.$$
Combining these two inequalities, we obtain
$$\norm{z_2} \leq \norm{(1+\lambda^2)z_2} \leq (1+|\lambda|)2\epsilon,$$
which implies
$$1=|||y_2||| = \norm{z_2}+\norm{\alpha z_2} \leq
(1+|\lambda|+\epsilon)(1+|\lambda|)2\epsilon,$$
a contradiction.

From this we deduce that $Y_1$ and $Y_2$ form a direct sum in $Y$.
Therefore $Y$ is not $\R$-HI. As $X$ is $HD_2$ as a real space
(\cite{F2}, Corollary 2) and
$Y \subset X$, it follows that $Y$ is $HD_2$ as a real space.
It follows that $Y$ is $\R$-quasimaximal in $X$, and therefore
$\C$-quasimaximal in $X$ (this follows immediately from the definition of
quasimaximality). As every
$\C$-linear subspace
$Y$ of
$X$ is
$\C$-quasimaximal in $X$, it follows that $X$ is H.I. as a complex
space.
\pff

\

This remark and the previous examples illustrate how various
the relations can be between real and complex structure in a complex H.I.
space. 

\

\paragraph{Question 34}
By Theorem \ref{valentin}, when $X$ is real H.I., there exists a division algebra $E$
isomorphic to
$\R$,
$\C$ or
$\HH$, such that for any $Y \subset X$, ${\cal L}(X)/{\cal S}(X)$ embeds into 
 ${\cal L}(Y,X)/{\cal S}(Y,X)$ which embeds into $E$.
It follows that for any $Y \subset X$,
$$\dim {\cal L}(X)/{\cal S}(X) \leq {\cal L}(Y,X)/{\cal S}(Y,X) \leq \dim
E.$$
 If $X=X_{GM}$ all these dimensions are equal to $1$. We provided 
examples $X(\C)$ and
$X(\HH)$ for which all these dimensions are equal to $2$ or to $4$
respectively. It remains open whether these dimensions
 may differ. For example, does there exist a real
H.I. Banach space $X$ such that every operator on $X$ is of the form
$\lambda Id_X +S$, but such that there exists an operator $T$ on some
subspace $Y$ of $X$ which is not of the form $\lambda i_{Y,X}+s$?
In this case, ${\cal L}(X)/{\cal S}(X)$ would be isomorphic to $\R$,
while $E$ would be complex or quaternionic.

\

Valentin Ferenczi,

Universit\'e Paris 6,

Equipe d'Analyse, Bo\^\i te 186,

4, place Jussieu, 75252 Paris Cedex 05,

France.

\

e-mail: ferenczi@ccr.jussieu.fr


\begin{thebibliography}{Wwww}

\bibitem{An} R. Anisca, {\em Subspaces of $L_p$ with more than one complex structure},
Proc. Amer. Math. Soc.  {\bf 131}  (2003),  no. 9, 2819--2829.

\bibitem{A} S. Argyros, {\em Ramsey
methods in analysis 
(Saturated and conditional structures in Banach spaces)},
 Advanced Courses in Mathematics,
CRM Barcelona. Birkhauser Verlag, Basel, (2005).

\bibitem{A2} S. Argyros, {\em A universal property of reflexive hereditarily
indecomposable Banach spaces}, Proc. A.M.S., {\bf 129} (2001), 3231--3239.  

 \bibitem{AM} S. Argyros and A. Manoussakis, 
{\em An indecomposable and unconditionally saturated Banach space},
Studia Math. {\bf 159} (2003), no. 1, 1--32.

\bibitem{AT} S. Argyros and A. Tolias, {\em Methods in the theory of
hereditarily indecomposable Banach spaces}, Mem. Amer. Math. Soc. {\bf
170} (2004), 806.

\bibitem{B} J. Bourgain, {\em Real isomorphic complex Banach
spaces need not be complex isomorphic},  Proc. Amer. Math. Soc. 
{\bf 96}  (1986),  no. 2, 221--226.

\bibitem{F1} V. Ferenczi, {\em Operators on subspaces of hereditarily
 indecomposable Banach spaces},
 Bull. London Math. Soc. {\bf 29} (1997), 338--344.


\bibitem{F2} V. Ferenczi, {\em Hereditarily finitely decomposable Banach spaces},
 Studia Mathematica {\bf 123} (2) (1997), 135--149.

\bibitem{F3} V. Ferenczi, {\em Quotient hereditarily indecomposable Banach spaces},  Canad. J. Math.  {\bf 51} 
(1999),  no. 3, 566--584.

\bibitem{GK} G. Godefroy and N. J. Kalton, {\em Lipschitz-free Banach spaces}, Studia Math.
{\bf 159} (2003), no.~1, 121--141. 



\bibitem{G0} W.T. Gowers, {\em A solution to Banach's hyperplane problem}, Bull.
London Math. Soc.  {\bf 26}  (1994),  no. 6, 523--530. 

\bibitem{G1} W.T. Gowers, {\em An infinite Ramsey theorem and some
Banach-space dichotomies}, Ann. of
Math. (2) {\bf 156} (2002), no. 3, 797--833.

\bibitem{GM1} W.T. Gowers and B. Maurey, {\em The unconditional basic
sequence problem}, J. Amer.
Math. Soc. {\bf 6} (1993),4, 851--874.

\bibitem{GM2} W.T. Gowers and B. Maurey, {\em Banach spaces with small spaces of
operators}, Math. Ann. {\bf 307} (1997), 543--568.

\bibitem{K} N. Kalton, {\em An elementary example of a Banach space not isomorphic to its
complex conjugate},  Canad. Math. Bull.  {\bf 38} (1995),  no. 2, 218--222.

\bibitem{LT} J. Lindenstrauss and L. Tzafriri, {\em Classical Banach
spaces}, Springer-Verlag, New York, Heidelberg, Berlin (1979).
  
\bibitem{M} B. Maurey, {\em Banach spaces with few operators}, Handbook of the geometry of
Banach spaces, vol. 2, edited by W.B. Johnson and J. Lindenstrauss, Elsevier, Amsterdam, 2002.

\bibitem{MU} S. Mazur and S. Ulam, {\em Sur les transformations
isom\'etriques d'espaces vectoriels norm\'es}, C. R. Acad. Sci.
Paris {\bf 194} (1932), 946--948.

\bibitem{S} S. Szarek, {\em On the existence and
uniqueness of complex structure and spaces with "few" operators}, 
Trans. Amer. Math. Soc.  {\bf 293}  (1986),  no. 1, 339--353.

\bibitem{S2} S. Szarek, {\em A superreflexive Banach space which does not admit complex structure},
  Proc. Amer. Math.
Soc.  {\bf 97}  (1986),  no. 3, 437--444.
\end{thebibliography}
\end{document}